\documentclass[reqno,11pt]{amsart}
\usepackage{a4wide}
\usepackage{amssymb}
\usepackage{amsmath}
\usepackage[pctex32]{graphics}
\usepackage{graphicx}

\newtheorem{theorem}{Theorem}
\newtheorem{corollary}[theorem]{Corollary}

\newtheorem{proposition}[theorem]{Proposition}

\theoremstyle{definition}
\newtheorem{definition}[theorem]{Definition}

\newtheorem{remark}[theorem]{Remark}

\newtheorem{example}[theorem]{Example}

\newcommand{\ds}{\displaystyle}

\newcommand{\be}{\begin{equation}}
\newcommand{\bd}{\begin{displaymath}}
\newcommand{\ed}{\end{displaymath}}
\newcommand{\ba}{\begin{eqnarray}}
\newcommand{\ea}{\end{eqnarray}}

\def\rr{\mathbb R}

\begin{document}

\title{Non ultracontractive heat kernel bounds by Lyapunov conditions}

\author{Fran\c cois Bolley}
\address{Ceremade, Umr Cnrs 7534, Universit\'e Paris-Dauphine, Place du Mar\'echal de Lattre de Tassigny, F-75775 Paris cedex 16} 
\email{bolley@ceremade.dauphine.fr}

\author{Arnaud Guillin}
\address{Institut Universitaire de France and Laboratoire de Math\'ematiques, Umr Cnrs 6620, Universit\'e Blaise Pascal, Avenue des Landais, F-63177 Aubi\`ere cedex}
\email{guillin@math.univ-bpclermont.fr}

\author{Xinyu Wang}
\address{Huazhong University of Science and Technology, Wuhan.}
\email{wang@math.univ-bpclermont.fr}

\date{\today}

\maketitle


\begin{abstract}
Nash and Sobolev inequalities are known to be equivalent to ultracontractive properties
of heat-like Markov semigroups, hence to uniform on-diagonal bounds on their kernel densities.  In non ultracontractive settings, such bounds can not hold, and (necessarily weaker, non uniform) bounds on the semigroups can be derived by means of weighted Nash (or super-Poincar\'e) inequalities.
The purpose of this note is to show how to check these weighted Nash inequalities in concrete examples of reversible diffusion Markov semigroups in $\rr^d$, in a very simple and general manner. We also deduce off-diagonal bounds for the Markov kernels of the semigroups, refining E.~B.~Davies' original argument. 
\end{abstract}

\bigskip

The (Gagliardo-Nirenberg-) Nash inequality in $\rr^d$ states that
\begin{equation}\label{nashrd}
\Vert f \Vert_{L^2(dx)}^{1+2/d} \leq C(d) \Vert \nabla f \Vert_{L^2(dx)} \, \Vert f \Vert_{L^1(dx)}^{2/d}
\end{equation}
for functions $f$ on $\rr^d$ and is a powerful tool when studying smoothing properties of parabolic partial differential equations on $\rr^d$. 

In a general way, let $(P_t)_{t \geq 0}$ be a symmetric Markov semigroup on a space $E$, with Dirichlet form $\mathcal E$ and (finite or not) invariant measure $\mu$. Then the Nash inequality
\begin{equation}\label{nashE}
\Vert f \Vert_{L^2(d\mu)}^{2+4/d} \leq \big[ C_1 \mathcal E (f,f) + C_2 \Vert f \Vert^2_{L^2(d\mu)} \big] \, \Vert f \Vert_{L^1(d\mu)}^{4/d}
\end{equation}
for a positive parameter $d$, or more generally
\begin{equation}\label{nashphi}
\Phi \left( \frac{\Vert f \Vert_{L^2(d\mu)}^{2}}{\Vert f \Vert_{L^1(d\mu)}^{2}} \right) \leq \frac{\mathcal E (f,f)}{\Vert f \Vert^2_{L^1(d\mu)}}
\end{equation}
for an increasing convex function $\Phi$, is equivalent,~up to constants and under adequate hypotheses on $\Phi$, to the ultracontractivity bound
$$
\Vert P_t f \Vert_{L^{\infty}(d\mu)} \leq n^{-1}(t) \, \Vert f \Vert_{L^1(d\mu)}, \quad  t > 0
$$
where $$\displaystyle n (t) = \int_t^{+ \infty} \frac1{\Phi(x)} dx;$$ for instance $n ^{-1} (t) \leq C t^{-d/2}$ for $0 < t \leq 1$ in the case of \eqref{nashE}. We refer in particular to~\cite{cks} and~\cite{davies} in this case, and  in the general case of \eqref{nashphi} to the seminal work~\cite{coulhon} by T.~Coulhon where the equivalence was first obtained; see also~\cite[Chap.~6]{bgl-book}. Let us observe that Nash inequalities are adapted to smoothing properties of the semigroup for small times, but can also be useful for large times. This is in turn equivalent to uniform~bounds on the kernel density of $P_t$ with respect to $\mu$, in the sense that for $\mu$-almost every $x$ in $E$ one can~write
\begin{equation}\label{ultra}
P_t f(x) = \int_E f(y) \, p_t(x,y) \, d\mu(y) \quad \textrm{with} \quad p_t(x,y) \leq  n^{-1}(t)
\end{equation}
for $\mu \otimes \mu$-almost every $(x,y)$ in $E \times E.$  Observe finally that \eqref{nashphi} is equivalent to its linearised~form
$$
\Vert f \Vert_{L^2(d\mu)}^{2} \leq u \, \mathcal E (f,f) + b(u) \, \Vert f\Vert_{L^1(d\mu)}^2, \qquad u<u_0
$$
for a decreasing positive function $b(u)$ related to $\Phi$: this form was introduced by F.-Y.~Wang~\cite{wang00} under the name of super-Poincar\'e inequality to characterize the generators $L$ with empty essential spectrum. For certain $b(u)$ it is equivalent to a logarithmic Sobolev inequality for~$\mu$, hence, to hypercontractivity only (and not ultracontractivity) of the semi\nobreak group.

\smallskip

Moreover, relevant Gaussian off-diagonal bounds on the density $p_t(x,y)$ for $x \neq y$, such~as
$$
p_t(x,y) \leq C(\varepsilon) \, t^{-d/2} \, e^{-d(x,y)^2/(4t (1+\varepsilon))}, \quad t>0
$$
for all $\varepsilon >0$, have first been obtained by E.~B.~Davies~\cite{davies-xy87} for the heat semigroup on a Riemannian manifold $E$, and by using a family of logarithmic Sobolev inequalities equivalent to \eqref{nashE}. Such bounds have been turned optimal in subsequent works, possibly allowing for $\varepsilon =0$ and the optimal numerical constant $C$ when starting from the optimal so-called entropy-energy inequality, and extended to more general situations: see for instance~\cite[Sect. 7.2]{bgl-book},~\cite{cks} and~\cite{davies} for a presentation of the strategy based on entropy-energy inequalities, and~\cite{boutayeb-coulhon-sikora}, \cite[Sect.~2]{coulhon-notes} and \cite{coulhon-sikora} and the references therein for a presentation of three other ways of deriving off-diagonal bounds from on-diagonal ones (namely based on an integrated maximum principle, finite propagation speed for the wave equation and a complex analysis argument).

\bigskip

In the more general setting where the semigroup is not ultracontractive, then the uniform bound~\eqref{ultra} cannot hold, but only (for instance on the diagonal) 
\begin{equation}\label{nonunif}
p_t(x,x) \leq  n^{-1} (t) V(x)^2
\end{equation}
for a nonnegative function $V$. Such a bound is interesting since it provides information on the semigroup : for instance if $V$ is in $L^2(\mu)$, then it ensures that $P_t$ is Hilbert-Schmidt, and in particular has a discrete spectrum. It has been shown to be equivalent to a {\it weighted} super-Poincar\'e inequality
\begin{equation}\label{SPIV}
\Vert f \Vert_{L^2(d\mu)}^{2} \leq u \, \mathcal E (f,f) + b(u) \, \Vert f V \Vert_{L^1(d\mu)}^2, \qquad u<u_0
\end{equation}
as in \cite{wang02}, where sharp estimates on high-order eigenvalues are derived, and, as in \cite{bbgm-nash}, to a {\it weighted} Nash inequality
\begin{equation}\label{nashphiV}
\Phi \left( \frac{\Vert f \Vert_{L^2(d\mu)}^{2}}{\Vert f V \Vert_{L^1(d\mu)}^{2}} \right) \leq \frac{\mathcal E (f,f)}{\Vert f V \Vert^2_{L^1(d\mu)}} \cdot
\end{equation}

\bigskip

The purpose of this note is twofold. First, to give simple and easy to check sufficient criteria on the generator of the semigroup for the weighted inequalities \eqref{SPIV}-\eqref{nashphiV} to hold : for this, we use Lyapunov conditions, which have revealed an efficient tool to diverse functional inequalities (see \cite{guillin-al-08-1} or \cite{cgww09} for instance) : we shall see how they allow to recover and extend examples considered in \cite{bbgm-nash} and \cite{wang02}, in a straightforward way (see Example \ref{exssexp}). Then, to derive off-diagonal bounds on the kernel density of the semigroup, which will necessarily be non uniform in our non ultracontractive setting. For this we refine Davies' original ideas of~\cite{davies-xy87}:  indeed, we combine his method with the (weighted) super-Poincar\'e inequalities derived in a first step, instead of the families of logarithmic Sobolev inequalities or entropy-energy inequalities used in the ultracontractive cases of~\cite{bgl-book},~\cite{cks} and~\cite{davies-xy87}-\cite{davies}; we shall see that the method recovers the optimal time dependence when written for (simpler) ultracontractive cases, and give new results in the non ultracontractive case (extending the scope of \cite{bbgm-nash} and~\cite{wang02}). Instead, we could have first derived on-diagonal bounds, such as~\eqref{nonunif}, and then use the general results mentionned above (in particular in~\cite{coulhon-sikora}) and giving off-diagonal bounds from on-diagonal bounds; but we will see here that, once the inequality \eqref{SPIV}-\eqref{nashphiV} has been derived, the off-diagonal bounds come without further assumptions nor much more effort than the on-diagonal ones.

\bigskip

To make this note as short and focused on the method as possible, we shall only present in detail the situation where $U$ is a $C^2$ function on $\rr^d$ with Hessian bounded by below, possibly by a negative constant, and such that $\int e^{-U} \, dx=1$. The differential operator $L$ defined by $L f = \Delta f - \nabla U \cdot \nabla f$ for $C^2$ functions $f$ on $\rr^d$ generates a Markov semigroup $(P_t)_{t \geq 0}$, defined for all $t \geq 0$ by our assumption on the Hessian of $U$. It is symmetric in $L^2(\mu)$ for the invariant measure $d\mu(x) = e^{-U(x)} \, dx.$ We refer to~\cite[Chap. 3]{bgl-book} for a detailed exposition of the background on Markov semigroups. Let us point out that the constants obtained in the statements do not depend on the lower bound of the Hessian of $U$, and that the method can be pursued in a more general setting, see Remarks~\ref{rmqvar},~\ref{rmsaut} and~\ref{rem11}.

\smallskip

We shall only seek upper bounds on the kernel, leaving lower bounds or bounds on the gradient aside (as done in the ultracontractive setting in~\cite[Sect. 2]{coulhon-notes},~\cite{davies} or~\cite{saloff} for instance). Let us finally observe that S.~Boutayeb, T.~Coulhon and A.~Sikora~\cite[Th.~1.2.1]{boutayeb-coulhon-sikora} have most recently devised a general abstract framework, including a functional inequality equivalent to the more general bound $p_t(x,x) \leq m(t,x)$ than~\eqref{nonunif}, and to the corresponding off-diagonal bound. The derivation of simple practical criteria on the generator ensuring the validity of such (possibly optimal) bounds is an interesting issue, that should be considered elsewhere.

\bigskip

\begin{definition}\label{defi} Let $\xi$ be a $C^1$ positive increasing function on $(0, + \infty)$ and $\phi$ be a continuous positive function on $\rr^d$, with $\phi (x) \to + \infty$ as $\vert x \vert$ goes to $+ \infty$. A $C^2$ map $W \geq 1$ on $\rr^d$ is a $\xi$-Lyapunov function with rate $\phi$ if there exist $b, r_0 \geq 0$ such that
for all $x \in \rr^d$
$$
\frac{L W}{\xi(W)} (x)\leq - \phi (x)+ b  \, 1_{\vert x \vert \leq r_0}.
$$
\end{definition}

We first state our general result:

\begin{proposition}\label{theogene}
In the notation of Definition \ref{defi}, assume that there exists a $\xi$-Lyapunov function $W$ with rate $\phi$. Then there exist $C$ and $s_0 >0$ such that for any positive continuous function $V$ on $\rr^d$
\begin{equation}\label{SPIgene}
\int f^2 \, d\mu \leq s \, \int \vert \nabla f \vert^2 \, (1\vee 1 / \xi ' (W)) \, d\mu \; + \; C \, g \circ \psi\Big( \frac{4}{s} \Big) \;  \Big( \frac{1}{s} \, \vee \, h \circ \psi \Big( \frac{4}{s} \Big) \Big)^{d/2} \, \Big( \int \vert f \vert \, V \, d\mu \Big)^2
\end{equation}
for all smooth $f$ and $s \leq s_0$. Here, for $r>0$
$$
g(r) = \sup_{ \vert x \vert \leq r} \frac{e^{U(x)}}{V(x)^2}, \quad h(r) = \sup_{ \vert x \vert \leq r} \vert \nabla U(x) \vert^2, \quad \psi(r) = \inf \{u \geq 0; \inf_{ \vert x \vert \geq u} \phi(x) \geq r\}.
$$
\end{proposition}

\smallskip

Here and below $a \vee b$ stands for $\max \{a, b \}.$

\begin{remark}\label{rmqvar}
 Proposition \ref{theogene} can be extended from the $\rr^d$ case to the case of a $d$-dimen\-sional connected complete Riemannian manifold $M$. If $M$ has a boundary $\partial M$ we then have to suppose that $\partial_n W \leq 0$ for the inward normal vector $n$ to the boundary, namely that the vector $\nabla W$ is outcoming at the boundary. We also have to check that a local super-Poincar\'e inequality holds: this inequality can easily be obtained in the $\rr^d$ case by perturbation of the Lebesgue measure, as in the proof below; for a general manifold, it holds if the injectivity radius of $M$ is positive or, with additional technical issues, if the Ricci curvature of $M$ is bounded below (see \cite{cgww09} or \cite{wang00} and the references therein).
\end{remark}

\medskip

\begin{corollary}\label{coro-weight}
Assume that there exist $c, \alpha$ and $\delta >0$ such that 
\begin{equation}\label{hypnablaU}
\nabla U(x) \cdot x \geq \frac{1}{c} \vert x \vert^{\alpha} -c, \qquad  \vert \nabla U(x) \vert \leq c \, (1 \vee \vert x \vert^{\delta})
\end{equation}
for all $x \in \rr^d$. Then for all $\gamma \geq 0, \gamma > 1 - \alpha$ and for $V(x) = (1+\vert x \vert^2)^{-\beta/2} e^{U(x)/2}$ with $\beta \in \rr$ there exist $C$ and $s_0 >0$ such that
\begin{equation}\label{SPIVomega}
\int f^2 \, d\mu \leq s \, \int \vert \nabla f \vert^2 \, (1\vee \vert x \vert^{2 \gamma}) \, d\mu + \frac{C}{s^p} \Big( \int \vert f \vert \, V \, d\mu \Big)^2
\end{equation}
for all smooth $f$ and $s \leq s_0$. Here
$$
p = \frac{2 \beta \vee 0 + d \, [\delta \vee (\alpha + \gamma-1)] }{2(\alpha + \gamma -1)}\cdot
$$
\end{corollary}

The first hypothesis in \eqref{hypnablaU} allows the computation of an explicit Lyapunov function $W$ as in Definition \ref{defi}, with an explicit map $\phi$, hence $\psi$ in~\eqref{SPIgene};  the second hypothesis is made here only to obtain an explicit map $h$ in~\eqref{SPIgene}, and then the $s^{-p}$ dependence as in the super Poincar\'e inequality~\eqref{SPIVomega}. 
Observe also from the proof that $s_0$ does not depend on $\beta$, and that the constant $C$ obtained by tracking in the proof its dependence on the diverse parameters would certainly be far from being optimal, as it is always the case in Lyapunov condition arguments.
\medskip

A variant of the argument leads to a (weighted if $\alpha <1$) Poincar\'e inequality for the measure~$\mu$, see Remark \ref{rmkpoincare} below. For $\alpha >1$ we can take $\gamma =0$ in Corollary~\ref{coro-weight}, obtaining a super-Poincar\'e inequality with the usual Dirichlet form and a weight $V$, and then 
non uniform off-diagonal bounds on the Markov kernel density of the associated semigroup:

\begin{theorem}\label{coro-noweight}
Assume that $\Delta U - \vert \nabla U \vert^2/2$ is bounded from above, and that there exist $c, \delta >0$ and $\alpha>1$ such that 
$$
\nabla U(x) \cdot x \geq \frac{1}{c} \vert x \vert^{\alpha} - c, \qquad \vert \nabla U(x) \vert \leq c \, (1 \vee \vert x \vert^{\delta})
$$
for all $x \in \rr^d$. 
Then for all $t>0$ the Markov kernel of the semigroup $(P_t)_{t \geq 0}$ admits a density $p_t(x,y)$ which satisfies the following bound : for all $\beta \in \rr$  and $\varepsilon>0$ there exists a constant~$C$ such~that
$$
 p_t(x,y) \leq  \frac{C}{t^{p}} \,  \frac{e^{U(x)/2} \, e^{U(y)/2}}{(1+\vert x \vert^2)^{\beta/2} (1+\vert y \vert^2)^{\beta/2}} \, e^{-\frac{\vert x-y \vert^2}{4(1+\varepsilon)t}}
$$
for all $0 < t \leq 1$ and almost every $(x,y)$ in $\rr^d \times \rr^d$. Here 
$$
p = \frac{2 \beta \vee 0 + d [ \delta \vee (\alpha -1)]}{2(\alpha -1)}\cdot
$$
\end{theorem}

 In particular, if $\delta \leq \alpha -1$ (as in Example \ref{exssexp} below), and for $\beta =0$, we obtain $p=d/2$ as in the ultracontractive case of the heat semigroup. The bound on the kernel density derived in Theorem~\ref{coro-noweight} will be proved for all $t>0$, but with an extra $e^{Ct}$ factor, hence is relevant only for small times. For larger times it can be completed as follows :

\begin{remark}\label{rempoincare}
In the assumptions and notation of Theorem~\ref{coro-noweight}, the measure $\mu$ satisfies a Poincar\'e inequality; there exists a constant $K>0$ and for all $t_0 >0$ there exists $C=C(t_0)$ such that for all $t\geq t_0$ and almost every $(x,y)$ in $\rr^d \times \rr^d$
$$
\vert p_t(x,y) - 1 \vert \leq C e^{-Kt} \frac{e^{U(x)/2} \, e^{U(y)/2}}{(1+\vert x \vert^2)^{\beta/2} (1+\vert y \vert^2)^{\beta/2}} \cdot
$$
\end{remark}

\begin{remark}\label{rmsaut}
Our method can be extended from the case of diffusion semigroups to more general cases. For instance the weighted Nash inequalities can be derived for discrete valued reversible pure jump process as then a local super-Poincar\'e inequality can easily be obtained. Then we should suppose that the Lyapunov condition holds for $\xi(w)=w$ and use \cite[Lem.~2.12]{cgww09} instead of \cite[Lem.~2.10]{cggr10} as in the proof below. Observe moreover that \cite{cks} has shown how to extend Davies' method \cite{davies-xy87} for off-diagonal bounds to non diffusive semigroups.
\end{remark}

\begin{example}\label{exssexp}
The measures with density $Z^{-1} e^{-u(x)^{\alpha}}$ for $\alpha > 0$ and $u$ $C^1$ convex such that $\int e^{-u} \, dx < + \infty$ satisfy the first hypothesis in \eqref{hypnablaU} in Corollary~\ref{coro-weight} and Theorem~\ref{coro-noweight}. Indeed, by \cite[Lem.~2.2]{guillin-al-08-1},
\begin{equation}\label{ucvx}
\nabla u(x) \cdot x \geq u(x) - u(0) \geq K \vert x \vert
\end{equation}
with $K>0$ and for $\vert x \vert$ large, so
$$
\nabla (u^{\alpha}) (x) \cdot x \, \vert x \vert^{-\alpha} = \alpha \, u(x)^{\alpha-1} \, \nabla u(x) \cdot x \vert x \vert^{-\alpha}
\geq \alpha \Big( \frac{u(x)}{\vert x \vert} \Big)^{\alpha} \Big( 1 - \frac{u(0)}{u(x)} \Big) \geq \alpha \Big( \frac{K}{2} \Big)^{\alpha} \, \frac{1}{2}>0
$$
for $\vert x \vert$ large.
 
 In particular for $u(x) = (1 + \vert x \vert^2)^{1/2}$, the hypotheses of Theorem \ref{coro-noweight} hold for $\alpha >1$ and $\delta = \alpha -1$ (and $U = u^{\alpha}$ has a curvature bounded by below), thus recovering the on-diagonal bounds given in \cite{wang02} (and \cite{bbgm-nash} in dimension $1$), and further giving the corresponding off-diagonal estimates. It was observed in \cite{bbgm-nash} that in the limit case $\alpha =1$ the spectrum of $-L$ does not only have a discrete part, so an on-diagonal bound such as $p_t(x,x) \leq C(t) V(x)^2$ with $V$ in $L^2(\mu)$ can not hold.  For $\alpha >2$ the semigroup is known to be ultracontractive (see~\cite{kkr} or~\cite[Sect. 7.7]{bgl-book} for instance), and adapting our method to this simpler case were $V=1$ leads to the corresponding off-diagonal  bounds, see Remark~\ref{remark-ultra}.  
\end{example}

\bigskip

For Cauchy type distributions, Proposition \ref{theogene} specifies as follows :

\begin{corollary}\label{coro-cauchy}
Assume that there exist $\alpha, c$ and $\delta>0$ such that 
$$
\liminf_{\vert x \vert \to + \infty} \nabla U(x) \cdot x \geq d + \alpha \qquad \textrm{and} \qquad \vert \nabla U(x) \vert \leq c \, (1 \vee \vert x \vert^{\delta}), \; x \in \rr^d.
$$
Then for all $\gamma > 2/\alpha$ and for $V(x) = (1+\vert x \vert^2)^{-\beta/2} e^{U(x)/2}$ with $\beta >0$ there exist $C$ and $s_0 >0$ such that
$$
\int f^2 \, d\mu \leq s \, \int \vert \nabla f \vert^2 \, (1\vee \vert x \vert^{2 +2/\gamma}) \, d\mu + \frac{C}{s^p} \Big( \int \vert f \vert \, V \, d\mu \Big)^2
$$
for all smooth $f$ and $s \leq s_0$. Here
$p = \gamma \, \beta + d \, [1 \vee (\delta \, \gamma)]/2.$
\end{corollary}

A variant of the argument also gives a weighted Poincar\'e inequality for the measure $\mu$, see Remark \ref{rmkpoincare}.

\begin{example}\label{excauchy}
The measures with density $Z^{-1} u(x)^{-(d+\alpha)}$ for $\alpha > 0$ and $u$ $C^2$ convex such that $\int e^{-u} \, dx < + \infty$ satisfy the first hypothesis in Corollary \ref{coro-cauchy}. Indeed, by \eqref{ucvx},
$$
\nabla U(x) \cdot x = (d+\alpha) \frac{\nabla u(x)}{u(x)} \cdot x \geq (d+\alpha) \frac{u(x) - u(0)}{u(x)} = (d+\alpha) \Big( 1- \frac{u(0)}{u(x)} \Big) \geq (d+\alpha) (1- \varepsilon)
$$
for any $\varepsilon>0$, and for $\vert x \vert$ large.  
In particular, choosing $u(x) = (1 + \vert x \vert^2)^{1/2}$, the hypotheses of Corollary \ref{coro-cauchy} hold with $\delta=0$ for the generalized Cauchy measures with density $Z^{-1}~(1+\vert x \vert^2)^{-(d+\alpha)/2}$, where $\alpha >0$.

\end{example}
\bigskip

\bigskip


The rest of this note is devoted to the proofs of these statements and further remarks.

\medskip

\noindent
{\bf Notation.} If $\nu$ is a Borel measure, $p \geq 1$ and $r>0$ we let $\Vert \cdot \Vert_{p, \nu}$ be the $L^p(\nu)$ norm and $\displaystyle \int_r f \, d\nu = \int_{\vert x \vert \leq r} f(x) \, d\nu(x).$ 
\bigskip


\noindent
{\bf Proof of Proposition \ref{theogene}.} We adapt the strategy of \cite[Th. 2.8]{cgww09}, writing
$$
\int \! f^2 \, d\mu = \int_r f^2 \, d\mu + \int_{ \vert x \vert > r} f^2 \, d\mu
$$
for $r>0$. By assumption on $\phi$ and $W$, and letting $\Phi(r) = \inf_{ \vert x \vert \geq r} \phi(x)$, the latter term is 
\begin{eqnarray*}
\int_{\vert x \vert > r} \! \! f^2 \, \phi \, \frac{1}{\phi} \, d\mu 
&\! \leq \!&
\frac{1}{\Phi( r)} \int_{ \vert x \vert > r} \! \! f^2 \, \phi \, d\mu \leq \frac{1}{\Phi( r)} \int f^2 \, \phi \,  d\mu
\\
&\! \leq \!&
\frac{1}{\Phi( r)} \left( \int \frac{-LW}{\xi(W)} f^2\, d\mu + b \int_{r_0} f^2 \, d\mu \right)
\leq \frac{1}{\Phi( r)} \left( \int \! \vert \nabla f \vert^2 \frac{1}{\xi'(W)} \, d\mu +  b \int_{r} \! f^2 \, d\mu \right)
\end{eqnarray*}
by integration by parts (as in \cite[Lem. 2.10]{cggr10}) and for $r \geq r_0$.

Hence for all $r \geq r_0$
\begin{equation}\label{est1}
\int f^2 \, d\mu
\leq
\Big(1+ \frac{b}{\Phi( r)} \Big)\int_r f^2 \, d\mu + \frac{1}{\Phi( r)} \int \vert \nabla f \vert^2 \omega \, d\mu
\end{equation}
for $\omega = 1\vee 1 / \xi'(W)$.

\bigskip

Now for all $r>0$ the Lebesgue measure on the centered ball of radius $r$ satisfies the following super-Poincar\'e inequality :
\begin{equation}\label{SPIball}
\int_r g^2\, dx \leq u \int_r \vert \nabla g \vert^2 \, dx + b(r,u) \left( \int_r \vert g \vert \, dx \right)^2
\end{equation}
for all $u>0$ and smooth $g$, where $b(r,u) = c_d(u^{-d/2} + r^{-d})$ for a constant $c_d$ depending only on $d$. For $r=1$ (say) this is a linearization of the (Gagliardo-Nirenberg-) Nash inequality
$$
\Big( \int_1 g^2 \, dx \Big)^{1+2/d} \leq c_d  \int_1 (g^2 + \vert \nabla g \vert^2) \, dx \; \;   \Big(\int_1 \vert g \vert \, dx \Big)^{4/d}
$$
for the Lebesgue measure on the unit ball; then the bound and the value of $b(r,u)$ for any radius~$r$ follow by homogeneity. 

Hence, for $g = f e^{-U/2}$,
\begin{eqnarray*}
\int_r f^2\, d\mu 
&=&
 \int_r (f e^{-U/2})^2 \, dx
 \\
&\leq&
2 \, u \int_r \vert \nabla f \vert^2 \, e^{-U} \, dx + \frac{u}{2} \int_r f^2 \vert \nabla U \vert^2 \, e^{-U} \, dx + b(r,u) \left( \int_r \frac{e^{U/2}}{V} \vert f \vert V e^{-U} dx \right)^2
\\
&\leq&
2 \, u \int \vert \nabla f \vert^2 \, \omega \, d\mu + \frac{u}{2} \, h( r) \int f^2 \, d\mu + b(r,u) \, g( r) \left( \int \vert f \vert V \, d\mu \right)^2.
\end{eqnarray*}

By \eqref{est1} and collecting all terms, it follows that for all $r\geq r_0$ and $u>0$
\begin{eqnarray}\label{est2}
\left( 1 - \frac{u}{2} h( r) \Big(1+ \frac{b}{\Phi( r)} \Big) \right) \int f^2 \, d\mu
&\leq&
 \left( \frac{1}{\Phi( r)} + 2u \Big(1+ \frac{b}{\Phi( r)} \Big) \right) 
\int \vert \nabla f \vert^2 \, \omega \, d\mu \nonumber
\\
&& + \Big(1+ \frac{b}{\Phi( r)} \Big) \, b(r,u) \, g( r) \left( \int \vert f \vert V \, d\mu \right)^2.
\end{eqnarray}

\smallskip

Now, for $r \geq r_0$, then $\Phi( r) \geq \Phi(r_0)$ so $1+ b / \Phi( r) \leq 1+ b / \Phi( r_0) := 1/k$. 
Hence, if $u \, h( r) \leq k$, then the coefficient on the left hand side in \eqref{est2} is greater than $1/2$; on the other hand, is $2 u \, \Phi( r) \leq k$, then the coefficient of the energy is smaller than $2/\Phi( r)$.

\medskip
Let now $s \leq s_0 := 4/\Phi(r_0)$ be given. Choosing $r := \psi(4/s)$ and then $u := k \,  \min \{1/ h( r), s/8 \}$, we observe that  $r \geq r_0$ and $u \leq ks_0/8$, so
$$
b(r,u) = c_d (u^{-d/2} + r^{-d}) \leq c_d(u^{-d/2} + r_0^{-d}) \leq C u^{-d/2}
$$ 
for a constant $C$ depending only on $d$ and $r_0$. Hence, for these $r$ and $u$, \eqref{est2} ensures the existence of constants $C$ and $s_0>0$ such that
$$
\int f^2 \, d\mu \leq s \int \vert \nabla f \vert^2 \, \omega \, d\mu + C g( r) \min \{1/h ( r), s \}^{-d/2}  \left( \int \vert f \vert V \, d\mu \right)^2
$$ 
for all weight $V$, $f$ and $s\leq s_0$. This concludes the proof of  Proposition \ref{theogene}. \hfill $\Box$

\bigskip

\noindent
{\bf Proof of Corollary \ref{coro-weight}.} Let $W \geq 1$ be a $C^2$ map on $\rr^d$ such that $W(x) = e^{a \vert x \vert^{\alpha}}$ for $\vert x \vert$ large, where $a>0$ is to be fixed later on. Then 
$$
L W (x) = a \alpha \, \vert x \vert^{2 \alpha - 2} W (x) \, \big[ (d + \alpha -2) \vert x \vert^{-\alpha} + a \alpha - \vert x \vert^{-\alpha} \, \nabla U(x) \cdot x \big]
$$
for $\vert x \vert$ large, by direct computation. Now
$$
a \alpha - \vert x \vert^{-\alpha} \, \nabla U(x) \cdot x \leq a \alpha - \frac{1}{2c} < 0
$$
for $\vert x \vert \geq (2c^2)^{1/\alpha}$ by the first assumption in \eqref{hypnablaU}, and for $a < (2c\alpha)^{-1}$, so for such an $a$ there exists a constant $C>0$ for which
\begin{equation}\label{LW}
LW (x) \leq - C \vert x \vert^{2\alpha - 2} W(x)
\end{equation}
for $\vert x \vert$ large.

\smallskip

In other words, if $\alpha >1$, then the Lyapunov condition of Definition \ref{defi} holds with $\xi(u) = u$ and $\phi(x) = C \vert x \vert^{2\alpha - 2}$: observe indeed  that $\phi(x) \to + \infty$ as $\vert x \vert \to + \infty$ for $\alpha >1$.

If the general case when possibly $\alpha \leq 1$, then we let $\gamma \geq 0$ and $\xi$ be a $C^1$ positive increasing map on $(0, + \infty)$ with $\xi(u) = u (\log u)^{-2\gamma/ \alpha}$ for, say, $u >e^{4\gamma/\alpha}$. Then, by \eqref{LW}, 
$$
\frac{LW}{\xi(W)} = \frac{LW}{a^{-2\gamma/ \alpha} \vert x \vert^{-2 \gamma} W} \leq - C \vert x \vert^{2(\alpha + \gamma -1)}
$$
for $\vert x \vert$ large, so the Lyapunov condition of Definition \ref{defi} holds with $\phi(x) = C \vert x \vert^{2(\alpha + \gamma-1)}$ if moreover $\gamma > 1 - \alpha$.

\smallskip

Hence Proposition \ref{theogene} ensures the super-Poincar\'e inequality \eqref{SPIgene} for $\mu$, with the weight $\omega (x) = 1 \vee 1/\xi'(W(x)) \leq C(1\vee \vert x \vert^{2 \gamma})$. We now choose $V(x) = (1+\vert x \vert^2)^{-\beta/2} e^{U(x)/2}$ with $\beta \in \rr$. Then, in the notation of  Proposition \ref{theogene},
$\psi(r) = C r^{1/2(\alpha+\gamma-1)}$ for all $r$,
$g( r) = (1+r^2)^{\beta} \leq 2^{\beta} r^{2 \beta}$ for $r \geq 1$ if $\beta >0$, and $g( r) \leq 1$ if $\beta \leq 0$; finally $h( r)  \leq c^2 r^{2\delta}$ under the second assumption in \eqref{hypnablaU}. This concludes the proof of Corollary \ref{coro-weight}. \hfill $\Box$


\bigskip

\noindent
{\bf Proof of Theorem \ref{coro-noweight}.}  It combines and adapts ideas from \cite{bbgm-nash} and \cite{davies-xy87}, replacing the family of logarithmic Sobolev inequalities used in \cite{davies-xy87} by the super-Poincar\'e inequalities~\eqref{SPIVomega}. The positive constants $C$ may differ from line to line, but will depend only on $U$ and $V$. 

The proof goes in the following several steps. 
\smallskip

{\bf 1.} Let $f$ be given in $L^2(\mu)$. With no loss of generality we assume that $f$ is non negative and $C^2$ with compact support, and satisfies $\displaystyle \int f \, V \, d\mu =1$. Let also $\rho>0$ be given and $\psi$ be a $C^2$ bounded map on $\rr^d$ such that $\vert \nabla \psi \vert \leq 1$ and $\vert \Delta \psi \vert \leq \rho$ (the formal argument would consist in letting $\psi(x) = x \cdot n$ for a unit vector $n$ in $\rr^d$). For a real number $a$ we also let $\displaystyle \varphi(x) = e^{a \psi(x)}$. We finally let $F(t) = \varphi^{-1} P_t(\varphi f)$.

\medskip

{\bf 1.1. Evolution of $\int F(t) V \, d\mu$.} We first observe that
\begin{equation}\label{ipp}
\frac{d}{dt} \int F(t) V \, d\mu = \int V \varphi^{-1} L P_t(\varphi f) \, d\mu = \int L (V \varphi^{-1}) P_t(\varphi f) \, d\mu
\end{equation}
by integration by parts. Indeed, following the proof of \cite[Cor. 3.1]{bbgm-nash}, two integrations by parts on the centered ball $B_r$ with radius $r>0$ ensure that
\begin{multline*}
\int_{B_r} V \varphi^{-1} \, L P_t (\varphi f) \, d\mu = \int_{B_r} L (V \varphi^{-1}) \, P_t (\varphi f) \, d\mu
\\
+ \int_{S^{d-1}} \Big[ P_t (\varphi f) (r \omega) \, \nabla (V \varphi^{-1}) (r \omega) \cdot n - (V \varphi^{-1}) (r \omega) \nabla P_t (\varphi f) (r \omega)  \cdot n \Big] e^{-U(r \omega)} r^{d-1} \, d\omega
\end{multline*}
for the inward unit normal vector $n$. Then a lower bound $\lambda \in \rr$ on the Hessian matrix of $U$ yields the commutation property and bound
$$
\vert \nabla P_t(\varphi f) \vert \leq e^{-\lambda t} P_t \vert \nabla (\varphi f) \vert \leq C e^{-\lambda t}
$$
by our assumption on $f$ and $\varphi$; moreover, on the sphere $\vert x \vert =r$, both
$(V \varphi^{-1}) (x) e^{-U(x)} \vert x \vert^{d-1}$ and
$\vert \nabla (V \varphi^{-1}) (x) \vert e^{-U(x)} \vert x \vert^{d-1}$ are uniformly bounded by $C (1+ r^C) \sup_{\vert x \vert = r} e^{-U(x)/2}.$ In turn this term goes to $0$ as $r$ goes to infinity since for instance
$$
U(x) - U(0) = \int_0^{r^{-1} (2c^2)^{1/\alpha}} \nabla U(tx) \cdot x \, dt + \int_{r^{-1} (2c^2)^{1/\alpha}}^1 \nabla U(tx) \cdot tx \, t^{-1} dt \geq \frac{r^{\alpha}}{C} - C
$$
for all $\vert x \vert =r$ large enough, by assumption on $ \nabla U$. Hence both boundary terms go to $0$ as $r$ goes to infinity; this proves \eqref{ipp}.

Then we observe that
$$
\frac{L (V \varphi^{-1})}{V \varphi^{-1}} = L ( \log (V \varphi^{-1})) + \vert \nabla \log (V \varphi^{-1}) \vert^2.
$$
But
$$
\log(V \varphi^{-1}) = \frac{U(x)}{2} - \beta \log < \!\!x\!\!> - a \psi(x)
$$
where $< \!\!x\!\!> = (1+\vert x \vert^2)^{1/2}$, so by direct computation
\begin{eqnarray*}
&& L ( \log (V \varphi^{-1})) + \vert \nabla \log (V \varphi^{-1}) \vert^2\\
&\!\!\!=&
\!\!\!\!
\frac{1}{2} \Delta U - \frac{1}{4} \vert \nabla U \vert^2  - \beta d \!< \!\!x\!\!>^{\!-2} \! + \beta(2 +\beta) \!< \!\!x\!\!>^{\!-4} \! \vert x \vert^2 - a \Delta \psi -2 a \beta < \!\!x\!\!>^{\!-2} \! \nabla \psi \cdot x+ a^2 \vert \nabla \psi \vert^2 \\
&\!\!\! \leq&
\!\!\!
C +  (\rho + \vert \beta \vert) \vert a \vert + a^2 := K
\end{eqnarray*}
for all $x$ if $\Delta U - \vert \nabla U \vert^2/2$ is bounded from above, $\vert \Delta \psi \vert \leq \rho$ and $\vert \nabla \psi \vert \leq 1$. Hence $L (V \varphi^{-1}) \leq K \, V \varphi^{-1}$, so
 $$
 \frac{d}{dt} \int F(t) V \, d\mu \leq K \int F(t) V \, d\mu,
 $$
 which implies
 $$
 \int F(t) V \, d\mu \leq e^{Kt} \int F(0) V \, d\mu = e^{Kt} \int f V \, d\mu = e^{Kt}.
 $$
 
 \medskip
 
 {\bf 1.2. Evolution of $y(t) := \int F(t)^2 \, d\mu$.} By integration by parts,
 \begin{eqnarray*}
 y'(t)
 &=&
 2 \int F \varphi^{-1} L P_t(\varphi f) \, d\mu 
 =
 -2 \int \nabla (F \varphi^{-1}) \cdot \nabla P_t(\varphi f) \, d\mu \\
 &=&
 -2 \int e^{-a \psi} (\nabla F - a F \nabla \psi)  \cdot e^{a \psi} (\nabla F + a F \nabla \psi) \, d\mu \\
 &=&
 -2 \int \vert \nabla F \vert^2 \, d\mu + 2 a^2 \int \vert \nabla \psi \vert^2 F^2 \, d\mu \leq -2 \int \vert \nabla F \vert^2 \, d\mu + 2 a^2 y(t)
 \end{eqnarray*}
 since $\vert \nabla \psi \vert \leq 1.$ But, for $\alpha >1$ we can take $\gamma =0$ in Corollary \ref{coro-weight}, so that
 $$
 \int \vert \nabla F \vert^2 \, d\mu \geq s^{-1} \Big[ \int F^2 \, d\mu - \frac{c}{s^p} \Big( \int F V d\mu \Big)^2 \Big]
 $$
 for any $s \leq s_0$ and a constant $c=c(U,V)$, and then
 $$
  \int \vert \nabla F \vert^2 \, d\mu \geq u \, y(t) - c \, u^{p+1} e^{2Kt}
 $$
for any $u(=s^{-1}) \geq u_0 := s_0^{-1}.$ Here $p = \big[2 \beta \vee 0  + d [ \delta \vee (\alpha -1)]\big]/ 2(\alpha -1).$ Hence
$$
y'(t) \leq - 2 (u \, y(t) - c \, u^{p+1} e^{2Kt}) + 2 a^2 y(t)
$$
for any $u \geq u_0$.
As long as $y(t) \geq c (p+1) u_0^p e^{2Kt}$ then the bracket  is minimal at 
$$
u = \Big(\frac{y(t)}{C(p+1)}\Big)^{1/p} e^{-2Kt/p} \geq u_0,
$$ 
and for this $u$ 
$$
y'(t) \leq -C e^{-2Kt/p} y(t)^{1+1/p} + 2 a^2 y(t).
$$
Then the map $z(t) = e^{-2a^2t} y(t)$ satisfies
$$
z'(t) \leq - C e^{-kt} z(t)^{1+1/p}
$$
where
$$
k:=\frac{2}{p} (K - a^2) = \frac{2}{p} (C + (\vert \beta \vert + \rho) \vert a \vert)) >0.
$$
It follows by integration that
$$
z(t)^{-1/p} \geq z(0)^{-1/p} + \frac{C}{p} \frac{1-e^{-kt}}{k} \geq \frac{C}{p} \frac{1-e^{-kt}}{k}
$$
so that
$$
y(t) = e^{2a^2t} z(t) \leq C e^{2a^2t} \Big(\frac{k}{1-e^{-kt}}\Big)^p \leq C e^{2a^2t}  \Big(\frac{e^{kt}}{t}\Big)^p = \frac{C}{t^p} e^{2Kt}.
$$
This last bound holds as long as $y(t) \geq c (p+1) u_0^p e^{2Kt}$, and then for all $t$ provided we take a possibly larger constant $C$ (still depending only on $U$ and $V$) in it and in the definition of $K$.

\medskip

{\bf 1.3.} In other words, for such a function $f$:
$$
\int \big(\varphi^{-1} P_t(\varphi f) \big)^2 \, d\mu \leq c(t)^2 
$$
for all $t>0$, where
$$
c(t)^2 = \frac{C}{t^p}e^{2Kt}, \quad K = C +  (\rho + \vert \beta \vert ) \vert a \vert + a^2.
$$

\medskip

{\bf 2. Duality argument.} Let $\varphi$ be defined as in step 1. By homogeneity and the bound $\vert P_t(\varphi f) \vert \leq P_t (\varphi \vert f \vert)$, it follows from step 1 that
$$
\int \big(\varphi^{-1} P_t(\varphi f) \big)^2 \, d\mu \leq c(t)^2 \Big( \int \vert f \vert V \, d\mu \Big)^2
$$
for all $t>0$ and all continuous function $f$ with compact support.

Let now $t>0$ be fixed, $Q$ defined by $Q f = \varphi^{-1} P_t(\varphi f)$ and $W = c(t) V$. Then
\begin{equation}\label{dua1}
\Vert Q f \Vert_{2, \mu} \leq \Vert f W \Vert_{1, \mu}.
\end{equation}
 The bound also holds for $-a$, so for $\varphi^{-1}$ instead of $\varphi$, so for $\varphi P_t(\varphi^{-1} f) = Q^*f$ where $Q^*$ is the dual of $Q$ in $L^2(\mu)$:
\begin{equation}\label{dua2}
\Vert Q^* f \Vert_{2, \mu} \leq \Vert f W \Vert_{1, \mu}.
\end{equation}
Following \cite[Prop. 2.1]{bbgm-nash} we consider $R f = W^{-1} Q (Wf)$, and its dual in $L^2(W^2 \mu)$, which is $R^* f =W^{-1} Q^*(Wf)$. Then it follows from \eqref{dua2} that
$$
\Vert R^* f \Vert_{2, W^2 \mu} \leq \Vert f  \Vert_{1, W^2\mu}
$$
and then by duality that
$$
\Vert R f \Vert_{\infty, W^2 \mu} \leq \Vert f  \Vert_{2, W^2\mu}.
$$
Moreover \eqref{dua1} ensures that
$$
\Vert R f \Vert_{2, W^2 \mu} \leq \Vert f  \Vert_{1, W^2\mu},
$$
so finally
$$
\Vert R^2 f \Vert_{\infty, W^2 \mu} \leq \Vert R f \Vert_{2, W^2 \mu} \leq \Vert f  \Vert_{1, W^2\mu}.
$$

As a consequence, by \cite[Prop. 1.2.4]{bgl-book} for instance, there exists a kernel density $r_2(x,y)$ on $\rr^d \times \rr^d$ such that $r_2(x,y)  \leq 1$ for $W^2\mu \otimes W^2 \mu$-almost every $(x,y)$ in $\rr^d \times \rr^d$ and
$$
R^2 f (x) = \int f (y) r_2(x,y) W(y)^2 \, d\mu(y)
$$
for all $f$ and $W^2\mu$-almost every $x$, hence (Lebesgue) almost every $x$ since $W \geq 1$ and $\mu$ has positive density. 

Observing that $P_t f = \varphi \, Q(\varphi^{-1}f)$ and $Q(g) = W \, R (W^{-1}g)$, it follows that 
\begin{eqnarray*}
P_{2t} f(x) 
&=&
\varphi (x) \, Q^2(\varphi^{-1}f) (x) = \varphi(x) W(x) R^2(W^{-1} \varphi^{-1}  f) (x)\\
&=&
\varphi(x) W(x) \int \varphi(y)^{-1} W(y) f(y) \, r_2(x,y) \, d\mu(y)
\end{eqnarray*}
for all $f$ and almost every $x$. Hence the Markov kernel of $P_{2t}$ has a density with respect to~$\mu$, given by
$$
p_{2t}(x,y) = \varphi(x) \varphi(y)^{-1} W(x) W(y) r_2(x,y) \leq e^{a(\psi(x)-\psi(y))} c(t)^2 V(x) V(y).
$$

\medskip

{\bf 3. Conclusion.} It follows from step 2 that the semigroup $(P_t)_{t \geq 0}$ at time $t>0$ admits a Markov kernel density $p_t(x,y)$ with respect to $\mu$, such that for all real number $a$ and all $C^2$ bounded map $\psi$ on $\rr^d$ with $\vert \nabla \psi \vert \leq 1$ and $\vert \Delta \psi \vert \leq \rho$, the bound
\begin{equation}\label{estpt}
p_{t}(x,y) \leq \frac{C}{t^p} V(x) V(y) \,  e^{Kt + a(\psi(x)-\psi(y))}
\end{equation}
hold for almost every $(x,y)$ in $\rr^d \times \rr^d$, with $K = C + (\rho + \vert \beta \vert) \vert a \vert + a^2.$

\smallskip

We now let $t>0$, $x$ and $y$ be fixed with $y \neq x$. Letting $r = \vert x \vert \vee \vert y \vert$ and $n = \frac{y-x}{\vert y-x \vert}$, we let $\psi$ be a $C^2$ bounded map on $\rr^d$ such that $\vert \nabla \psi \vert \leq 1$ and $\vert \Delta \psi \vert \leq \rho$ everywhere, and such that $\psi(z) = z \cdot n$ if $\vert z \vert \leq r$. For instance we let $h(z)$ be a $C^2$ map on $\rr$ with $h(z) = z$ if $\vert z \vert \leq r$, $h$ constant for $\vert z \vert \geq R$ and satisfying $\vert h' \vert \leq 1$ and $\vert h''\vert \leq \rho$, which is possible for $R$ large enough compared to $\rho^{-1}$; then we let $\psi(z) = h(z \cdot n)$.

Such a map $\psi$ satisfies $\psi(x) - \psi(y) = - \vert x-y \vert$, so the quantity in the exponential in~\eqref{estpt}~is
$$
-a \vert x-y \vert + (a^2+ (\rho + \vert \beta \vert) \vert a \vert) t + C t.
$$
Since $\rho >0$ is arbitrary (and the constant $C$ depends only on $U$ and $V$) we can let $\rho$ tend to~$0$. Then we use the bound
$$
\vert a \beta \vert  \leq \varepsilon a^2 + \frac{1}{4 \varepsilon} \beta^2
$$
and optimise the obtained quantity by choosing $a=\vert x-y \vert/(2t(1+\varepsilon))$, leading to the bound
$$
p_{t}(x,y) \leq \frac{C}{t^p} V(x) V(y) \, \exp \Big[ - \frac{\vert x-y \vert^2}{4(1+\varepsilon t)} + \frac{\beta^2 t}{4 \varepsilon} + C t \Big]
$$
and concluding the argument. \hfill $\Box$

\bigskip

\begin{remark}\label{rem11}
The computation in steps 1 and 2 of this proof could be written in the more general setting of a reversible diffusion generator $L$ on a space $E$ with carr\'e du champ $\Gamma$, under the assumptions $\Gamma(\psi) \leq 1, L (V \varphi^{-1}) \leq K V \varphi^{-1}$ and
$$
\int g^2 \, d\mu \leq s \int \Gamma(g) \, d\mu  + C s^{-p} \Big( \int \vert g \vert \, V \, d\mu \Big)^2
$$
for all $g$ and $s \leq s_0$. Then step 3 would yield the corresponding bound with $\vert x-y \vert$ replaced by the intrinsic distance
$$
\rho(x,y) = \sup \big\{ \vert \psi(x) - \psi(y) \vert; \Gamma(\psi) \leq 1 \big\}.
$$
\end{remark}

\bigskip


\noindent
{\bf Proof of Remark \ref{rempoincare}.}
First observe that, under the first hypothesis in \eqref{hypnablaU} in Corollary~\ref{coro-weight}, with $\alpha \geq 1$, then \eqref{LW} in the proof of this corollary ensures that
$$
\frac{LW}{W} (x) \leq - C + b \, 1_{\vert x \vert \leq r_0}
$$
for all $x$ and for positive constants $b$ and $C$. This is a sufficient Lyapunov condition for $\mu$ to satisfy a Poincar\'e inequality (see \cite[Th. 1.4]{guillin-al-08-1}).

Then we can adapt an argument in \cite[Sect. 7.4]{bgl-book}, that we recall for convenience. We slightly modify the notation of the proof of Corollary~\ref{coro-weight}, letting $a=0$ (hence $\varphi =1$), and for any $t \geq 0$, letting $R_t f= V^{-1} \, P_t (V f)$ and $R_{\infty} f = V^{-1} \, \ds \int  V f \, d\mu$. Then, by the Poincar\'e inequality for $\mu$, there exists a constant $K>0$ such that
\begin{eqnarray*}
\Vert R_{t} f - R_{\infty} f \Vert^2_{2, V^2 \mu} = \int \Big\vert P_t(Vf) - \int V f \, d\mu \Big\vert^2 \, d\mu 
&\leq&
 e^{-2Kt} \int \Big\vert Vf - \int V f \, d\mu \Big\vert^2 \, d\mu 
 \\
 &\leq&
  e^{-2Kt} \int \vert Vf \vert^2 \, d\mu
= e^{-2Kt} \Vert f \Vert^2_{2, V^2 \mu}
\end{eqnarray*}
for all $t \geq 0$. Moreover
$$
\Vert R_{t_0} f \Vert_{2, V^2 \mu} \leq c(t_0) \, \Vert f \Vert_{1, V^2 \mu}, \quad \Vert R_{t_0} f \Vert_{\infty, V^2 \mu} \leq c(t_0) \, \Vert f \Vert_{2, V^2 \mu}
$$
for all $t_0 >0$, by step 2 in the proof of Corollary~\ref{coro-weight}, so
\begin{multline*}
\Vert R_{t+2 t_0} f - R_{\infty} f \Vert_{\infty, W^2 \mu} 
=
 \Vert R_{t_0} R_t R_{t_0} f - R_{t_0} R_{\infty} R_{t_0} f \Vert_{\infty, W^2 \mu} 
\\
\leq
 c(t_0) \Vert R_t R_{t_0} f - R_{\infty} R_{t_0} f \Vert_{2, W^2 \mu}
\leq c(t_0) e^{-Kt} \Vert R_{t_0} f \Vert_{2, W^2 \mu}
\leq c(t_0)^2 e^{-Kt} \Vert f \Vert_{1, W^2 \mu}
\end{multline*}
for all $t\geq 0$ and $t_0>0$. Changing $t+ 2 t_0$ into $t \geq  t_0$ and writing the kernel density of $R_{t}  - R_{\infty}$ in terms of the kernel density $p_t(x,y)$ of $P_t$ lead to the announced bound on the density $p_t(x,y) - 1$. \hfill $\Box$
\bigskip

\begin{remark}\label{remark-ultra}
The proof of  Theorem~\ref{coro-noweight} simplifies in ultracontractive situations where one takes $1$ as a weight $V$. For instance, for the heat semigroup on $\rr^d$, one starts from the (non weighted) super Poincar\'e inequality for the Lebesgue measure on $\rr^d$:
$$
\int F^2 dx \leq u \int \vert \nabla F \vert^2 dx + c_d u^{-d/2} \Big( \int \vert F \vert dx \Big)^2
$$
for all $u>0$ and for a constant $c_d$ depending only on $d$ (a linearization of the Nash inequality~\eqref{nashrd} for the Lebesgue measure on $\rr^d$, which can also be recovered by letting $r$ go to $+ \infty$ in~\eqref{SPIball}). The constant $K$ obtained in step 1.1 is $K = \rho \vert a \vert + a^2$, and the very same argument leads to the (optimal in $t$) off-diagonal bound
$$
p_t(x,y) \leq C\, t^{-d/2} e^{-\frac{\vert x-y \vert^2}{4t}}
$$
for all $t>0$ (also derived in \cite[Sect. 7.2]{bgl-book} for instance with the optimal $C = (4 \pi)^{-d/2}$ when starting from by the Euclidean logarithmic Sobolev inequality).

The argument can also be written in the ultracontractive case of $U(x) = (1+ \vert x \vert^2)^{\alpha/2}$ with $\alpha >2$ : one starts from the super Poincar\'e inequality
$$
\int F^2 d\mu \leq u \int \vert \nabla F \vert^2 d\mu + e^{c\big(1+u^{-\frac{\alpha}{2\alpha-2}} \big)} \Big( \int \vert F \vert d\mu \Big)^2
$$
for all $u>0$ (see \cite[Cor.~2.5]{wang00} for instance). Then one obtains the off-diagonal bound
$$
p_t(x,y) \leq e^{C \big(1+t^{-\frac{\alpha}{\alpha-2}} \big)} \, e^{-\frac{\vert x-y \vert^2}{4t}}
$$
for all $t>0$, also derived in \cite[Sect. 7.3]{bgl-book} by means of an adapted entropy-energy inequality. 
\end{remark}


\bigskip

\noindent
{\bf Proof of Corollary \ref{coro-cauchy}.} Following the strategy of the proof of Corollary \ref{coro-weight}, we observe that $LW(x) \leq -C \vert x \vert^{a-2}$ with $C>0$ for $\vert x \vert$ large, if $W(x) = \vert x \vert^{a}$ with $a < 2 +\alpha$. Then, for $\xi(r) = r^{1-b}$, $LW/\xi(W) \leq - C \vert x \vert^{ab-2}$, and we are in the framework of  Proposition~\ref{theogene} provided $ab>2, a < 2+\alpha$ and $0 \leq b < 1$. We finally let $\gamma = 2/(ab-2).$ \hfill $\Box$

\bigskip

\begin{remark}\label{rmkpoincare}
We have seen in Remark \ref{rempoincare} that, by \eqref{LW}, the sole first hypothesis in \eqref{hypnablaU} in {\bf Theorem} \ref{coro-noweight}, for $\alpha \geq 1$, ensures a Poincar\'e inequality for $\mu$.

In the case when $0 < \alpha < 1$, let us take $\gamma = 1 - \alpha$ and replace $\omega(x) = C (1 \vee \vert x \vert^{2 \gamma})$ by a $C^1$ map $\omega(x)$, still equal to $C \vert x \vert^{2 \gamma}$ for $\vert x \vert $ large, and let 
$$
L^{\omega} f = \omega \, \Delta f - (\omega \nabla U - \nabla \omega) \cdot \nabla f = \omega L f + \nabla \omega \cdot \nabla f
$$
be the generator with symmetric invariant measure $\mu$. Then, still under the first hypothesis in \eqref{hypnablaU}, \eqref{LW} implies that
\begin{equation}\label{lomega}
\frac{L^{\omega} W}{W} = \omega \frac{LW}{W} + \nabla \omega \cdot \nabla \log W \leq  -C + c \vert x \vert^{-\alpha} \leq -C'
\end{equation}
for $\vert x \vert$ large and $C'>0$. This now leads to a Poincar\'e inequality for $\mu$ and for the energy
$$
-\int f \, L^{\omega} f \, d\mu = \int \vert \nabla f \vert^2 \, \omega \, d\mu \leq C \int  \vert \nabla f \vert^2 (1+\vert x \vert^{2(1- \alpha)}) \, d\mu,
$$
hence deriving, under the sole first hypothesis in \eqref{hypnablaU}, the weighted Poincar\'e inequality obtained in \cite[Prop.~3.6]{cggr10} for measures as in Example \ref{exssexp}. 

\smallskip

Likewise, under the sole first hypothesis in Corollary \ref{coro-cauchy}, the weight $\omega (x) = 1+ \vert x \vert^2$ satisfies the conclusion of \eqref{lomega}, leading to a Poincar\'e inequality for $\mu$ and for the energy
$$
-\int f \, L^{\omega} f \, d\mu = C \int \vert \nabla f \vert^2 (1+\vert x \vert^2) \, d\mu,
$$
 that is, to the weighted Poincar\'e inequality 
$$
\int \Big( f - \int f \, d\mu \Big)^2 d\mu \leq C \int \vert \nabla f \vert^2 (1+\vert x \vert^2) \, d\mu
$$
for measures as in Example \ref{excauchy}, with $\alpha >0$. For measures $\mu$ with density $Z^{-1} (1 + \vert x \vert^2)^{-(d + \alpha)/2}$, this inequality has been derived in~\cite[Prop.~3.2]{cggr10} for any $\alpha >0$. The dependence in $d$ and $\alpha$ of the constants $C$ obtained here and in~\cite{cggr10} is not tractable. An explicit constant $C = C(d, \alpha)$ has been derived in \cite[Th.~3.1]{bobkovledoux-cauchy} for $\alpha \geq d$, and improved in a sharp manner in \cite[Cor. 14]{nguyen} for $\alpha \geq d+2$.
 \end{remark}


\bigskip

\footnotesize \noindent\textit{Acknowledgments.} The authors warmly thank I. Gentil for enlightening discussion on this and related questions. They are grateful to the two referees for a careful reading of the manuscript and most relevant
comments and questions which helped improve the presentation of the paper. This note was partly
written while they were visiting the Chinese Academy of Sciences in Beijing; it is a pleasure for them to thank this
institution for its kind hospitality. The first two authors also acknowledge support from the French ANR-12-BS01-0019 STAB project.

\

\bibliographystyle{plain}

\end{document}